\documentclass[10pt]{conm-p-l}

\usepackage{amsmath,amsthm,amsfonts,amssymb,amscd}
\usepackage{pifont,pst-all}
\usepackage{lineno,color}
\usepackage{flafter}
\usepackage{graphicx}

\newtheorem{teo}{Theorem}

\newtheorem{pro}{Proposition}

\theoremstyle{definition}
\newtheorem{defi}{Definition}
\theoremstyle{remark}
\newtheorem{rem}{Remark}
\theoremstyle{definition}

\theoremstyle{definition}

\def    \bfC    {{\mathbf C}}

\def    \bfone  {{\mathbf 1}}
\def    \ssminus        {{\smallsetminus}}

\newcommand{\Q}{\mathcal Q}         

\newcommand{\dsum}{\displaystyle\sum}  
\newcommand{\dint}{\displaystyle\int}  

\newcommand{\Z}{\mathbb{Z}}       
\newcommand{\R}{\mathbb{R}}       
\newcommand{\C}{\mathbb{C}}       
\newcommand{\lbdle}{\mathbb{L}}   
\newcommand{\lclass}{\mathbf{L}}  
\newcommand{\todd}{\mathbf{Td}}   
\newcommand{\Qy}{\mathbf{Q}}  
\newcommand{\ind}{\mathrm{ind}}   

\newcommand{\e}{\mathrm{e}}       
\newcommand{\T}{\mathrm{T}}       
\newcommand{\CP}{\C\mathrm{P}}    

\newcommand{\liet}{\mathfrak{t}}     

\begin{document}

\bibliographystyle{hamsplain}

\title[A weighted version of quantization commutes with reduction]
{A weighted version of quantization commutes with reduction for a toric manifold}

\author[Jos\'e Agapito]{Jos\'e Agapito}

\address{Department of Mathematics, UC Santa Cruz, Santa Cruz, CA 95064, USA.}
\email{jarpepe@math.ucsc.edu}
\thanks{This work was partially supported by a UC Regents Fellowship the author
received during Spring 2003.}

\begin{abstract}
We compute explicitly the equivariant Hirzebruch
$\chi_y$-characteristic of an equivariant complex line bundle over
a toric manifold and state a weighted version of the quantization
commutes with reduction principle in symplectic geometry. Then, we
give a weighted decomposition formula for any simple polytope in
$\R^n$. This formula generalizes a polytope decomposition due to
Lawrence \cite{Law} and Varchenko \cite{Var} and extends a
previous weighted version obtained by Karshon, Sternberg and
Weitsman \cite{KSW1}.
\end{abstract}
\maketitle

\section{Introduction}
\label{se:intro}
A convex polytope is a convex bounded polyhedron in $\R^n$. We say
that a polytope is simple if there are exactly $n$ edges emanating
from each of its vertices. In addition, if the lines along the
edges of a simple polytope are generated by a $\Z$-basis of the
lattice $\Z^n$, the polytope is called regular.

Let $M$ be a toric manifold and let $\T$ be the torus acting on
$M$. The image of $M$ under the moment map of the torus action is
a regular polytope $\Delta$ which lies on $\liet^*$ (the dual of
the Lie algebra $\liet$ of $\T$). Each integral lattice point in
$\liet^*$ corresponds to a character of a representation of $\T$.
Counting lattice points in $\Delta$ is closely related to the
quantization commutes with reduction principle. For the Dolbeault
quantization, this principle states that each character $\alpha$
of $\T$ occurs in the quantization of $M$ with multiplicity 1 or 0
depending on whether the lattice point representing $\alpha$ is in
$\Delta$ or not. (c.f. \cite{Gin/Gui/Kar}, \cite{Gui}.)\par

We use the Atiyah-Bott Berline-Vergne localization formula in
equivariant cohomology and a notion of polarization to give a
weighted version of the quantization commutes with reduction
principle for a toric manifold. This is expressed as a way of
counting lattice points in $\Delta$ with \emph{weights}. From a
combinatorial point of view, we can count lattice points in
$\Delta$ using a polytope decomposition due to Lawrence \cite{Law}
and Varchenko \cite{Var}. (See Figure \ref{ordi-pic} in Subsection
\S\S\ref{subse:simplepolytopes}.) We give a weighted decomposition
formula for any simple polytope that generalizes the Lawrence and
Varchenko polytope decomposition. This formula also extends a
previous weighted version obtained by Karshon, Sternberg and
Weitsman. (See \cite{KSW2} and \cite{KSW1}.)

\subsection*{Acknowledgements}
I want to thank Jonathan Weitsman for his illuminating
explanations and comments throughout this article. I would like to
thank Shlomo Sternberg for his helpful comments on this paper.
Finally, I would like to thank the referee for several useful
suggestions.\par

\section{Definitions and Facts}
\label{se:definitions-facts}
For the reader's convenience, we collect in this section some
definitions, facts and remarks to be used later on in this
paper.\par

\subsection{Ordinary polar decomposition}
\label{subse:simplepolytopes} A convex \textbf{polytope} $\Delta$
in $\R^n$ is a \textbf{convex bounded polyhedron} given by an
intersection of half-spaces $\Delta=\cap H_i$, where
$$H_i=\{x\in\R^n\,|\,\langle x,u_i\rangle\ge\lambda_i\, ,\lambda_i\in\R\}
\qquad\mbox{for } i=1,\cdots,N>n,$$ with $u_i\in
(\R^n)^*\simeq\R^n$ and where $N$ denotes the number of
$(n-1)$-dimensional faces (the facets) of $\Delta$. The vector
$u_i$ can be thought of as the inward normal to the $i$-th facet
of $\Delta$ contained in the hyperplane $\langle
x,u_i\rangle=\lambda_i$. We assume that $\Delta$ is obtained with
the smallest possible $N$. We also assume that our polytopes are
always compact and convex.\par

A polytope is \textbf{simple} if there are exactly $n$ edges emanating from each vertex and is
\textbf{regular\footnote{our regular polytopes are also known as smooth or Delzant polytopes
\cite{Del}, \cite{Gui}.}} if, in addition, the edges emanating from each vertex lie along lines
which are generated by an integral basis of $\Z^n$. If the vertices of a polytope are in the
lattice $\Z^n$, we call it \textbf{integral}.\par

Let $\Delta$ be a simple polytope in $\R^n$. Then, we have $n$
edges emanating from each of its vertices. Let $v$ be a vertex. We
choose $\alpha_{1,v},\ldots,\alpha_{n,v}$ in the direction of the
edges pointing from $v$ to other vertices of $\Delta$. Up to a
positive scalar and permutation, these \emph{edge vectors} are
uniquely determined. When $\Delta$ is regular and integral, we can
choose $\alpha_{1,v},\ldots,\alpha_{n,v}$ to have integer entries.
In this case, we shall impose the additional normalization
condition that the edge vectors be primitive elements of the
lattice $\Z^n$; namely, none of them can be expressed as a
multiple of a lattice element by an integer greater than one. This
fixes our choice of edge vectors for such polytopes.\par

The tangent cone $\bfC_{v}$ at $v$ is
$$\bfC_{v}=\{v+r(x-v)\,|\, r\ge 0, x\in\Delta\}=v+\dsum_{j=1}^n
\R_{\ge 0}\alpha_{j,v}.$$

Notice that we can always find a vector $\xi$ such that
$\langle\alpha_{j,v},\xi\rangle\neq 0$ for all $v$ and $j$. To see
this, we can think of an edge vector as the normal vector of a
hyperplane in $\R^n$ passing through the origin. Since there are a
finite number of edge vectors, we have a finite number of such
hyperplanes. They split $\R^n$ into several connected components.
We pick $\xi$ such that it is not perpendicular to the boundary of
these components.\par

\begin{defi} A \textbf{polarizing vector} $\xi$ for a simple
polytope $\Delta\subset\R^n$ is a vector in $\R^n$ such that
$\langle\alpha_{j,v},\xi \rangle\neq 0$ for all vertices $v$ and
all edge vectors $\alpha_{1,v}, \ldots,\alpha_{n,v}$.
\end{defi}

At each vertex $v$, we define the polarized vectors
$\alpha^\sharp_{j,v}$ by
\begin{equation}\label{polarization}
\alpha^\sharp_{j,v}=\left\{%
\begin{array}{rrc}
  \alpha_{j,v} & \mbox{if }\langle\alpha_{j,v},\xi\rangle<0 & \quad\mbox{(unflipped)} \\
  -\alpha_{j,v} & \mbox{if }\langle\alpha_{j,v},\xi\rangle>0 & \quad\mbox{(flipped)} \\
\end{array}%
\right.\,.
\end{equation}
This process is called \textbf{polarization}. Then, the polarized
tangent cone $\bfC^\sharp_{v}$ at $v$ is
$$\bfC^\sharp_{v}=v+\dsum_{j=1}^n \R_{\ge 0}\alpha^\sharp_{j,v}.$$

We have the following (adapted) decomposition of a simple polytope
due to Lawrence \cite{Law} and Varchenko \cite{Var}.

\begin{pro}\label{ord-decomp} Let
$\Delta$ be a simple polytope in $\R^n$. For any choice of
polarizing vector, we have
\begin{equation}\label{oldformula}
\bfone_{\Delta}=\sum_{v}(-1)^{\#v}\bfone_{\bfC^\sharp_{v}},
\end{equation}
where the sum is over all vertices $v$ of the polytope and the
symbols $\bfone_{\Delta}$ and $\bfone_{\bfC^\sharp_v}$ denote the
ordinary characteristic functions over $\Delta$ and over the
polarized tangent cones $\bfC^\sharp_v$ (possibly with some facets
removed). In this formula, $\# v$ denotes the number of edge
vectors at $v$ flipped by the polarization process
\eqref{polarization}.
\end{pro}

\begin{figure}[h]
\centering
\includegraphics[scale=.85]{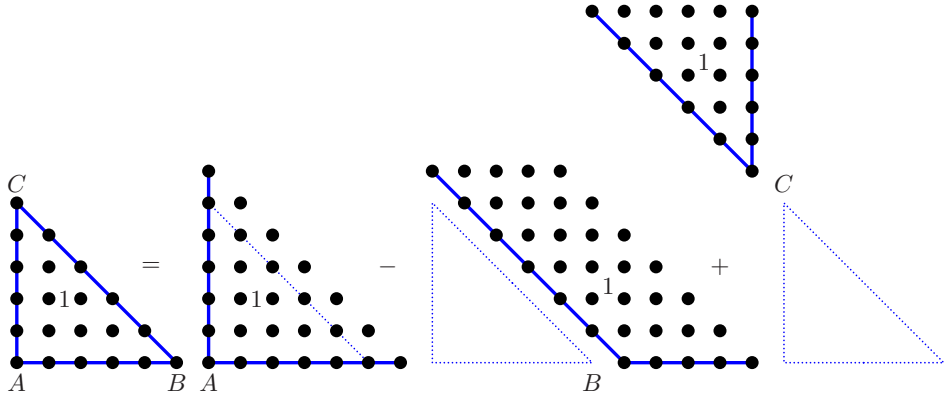}
\caption{Ordinary polar decomposition of a lattice
triangle}\label{ordi-pic}
\end{figure}

Lawrence \cite{Law} and Varchenko \cite{Var} formulas actually
hold for any simple convex polyhedron of full dimension in $\R^n$.
In the case of lattice polytopes, $\Delta\cap\Z^n$, we avoid
removing certain facets from the cones $\bfC^\sharp_v$, by
\emph{shifting} some of the cones. We illustrate this
decomposition in Figure \ref{ordi-pic}.\par

\subsection{Some classical formal power series}\label{subse:powerseries}
We will pay attention to three formal power series.\par

The Todd function $\todd(x)$ is defined as
$$\todd(x)=\dfrac{x}{1-\e^{-x}}=1+\dfrac{1}{2}x+\dfrac{1}{12}x^2-
\dfrac{1}{720}x^4+\cdots\,\,.$$ Since
$$\todd(-x)=\dfrac{x}{\e^{x}-1}=\e^{-x}\todd(x),$$
we can formally average $\todd(x)$ and $\todd(-x)$ to get
$$\begin{array}{lcl}
\dfrac{1}{2}(\todd(x)+\todd(-x))&=&\dfrac{1}{2}\todd(x)(1+\e^{-x})%
=(x/2)\dfrac{1+\e^{-x}}{1-\e^{-x}}\\[2ex]
&=&(x/2)\dfrac{e^{x/2}+\e^{-x/2}}{\e^{x/2}-\e^{-x/2}}=\dfrac{x/2}{\tanh(x/2)},
\end{array}$$

\noindent which is known as the $\hat{\lclass}(x)$ function. Both
$\todd(x)$ and $\hat{\lclass}(x)$ are particular instances of the
more general formal power series
$$Q(y,x)=\dfrac{x(1+y)}{1-\e^{-x(1+y)}}-yx=\dfrac{x(1+y\e^{-x(1+y)})}
{1-\e^{-x(1+y)}}.$$ This is the Hirzebruch function. When $y=0$ we
get $\todd(x)$ and if $y=1$ and $x$ is replaced by $x/2$, we get
$\hat{\lclass}(x)$. The Hirzebruch function can also be written as
$$Q(y,x)=\dfrac{1}{1+y}\todd(x(1+y))+\dfrac{y}{1+y}\todd(-x(1+y)).$$

\noindent We replace $x$ by $x/(1+y)$ above and denote the
resulting formula by $\Qy_y(x)$. We have a weighted average of
$\todd(x)$ and $\todd(-x)$. For any $y\neq -1$ it follows that
\begin{equation*}
\Qy_y(x)=\dfrac{1}{1+y}\todd(x)(1+y\e^{-x}).
\end{equation*}

  These classical formal power series define well known topological invariants;
namely, characteristic classes associated to vector bundles over a
manifold. Let $M$ be a compact manifold and let $E\to M$ be a
complex vector bundle over $M$. By the Splitting principle in
topology (c.f. \cite{BT}), we can always assume that $E$ splits
into a direct sum of complex line bundles,
$$E=L_1\oplus\cdots\oplus L_r.$$

Then the Todd class $\todd(E)$ of $E$ is given by
\begin{equation}\label{todd}
\todd(E)=\prod_{k=1}^r\dfrac{c_1(L_k)}{1-\e^{-c_1(L_k)}},
\end{equation}

\noindent where $c_1(L_1),\ldots,c_1(L_r)$ are the first Chern
classes of the line bundles $L_1,\ldots,L_r$.

\noindent Similarly, the $\hat{\lclass}$-class of $E$ and more
generally the Hirzebruch $\Qy_y$-class of $E$ are respectively
defined by
\begin{equation}\label{ordi-classes}
\begin{split}
&\hat{\lclass}(E)=\todd(E)\prod_{k=1}^r\frac{1}{2}(1+\e^{-c_1(L_k)})
\quad\mbox{and}\\
&\Qy_y(E)=\todd(E)\prod_{k=1}^r\frac{1}{1+y}(1+y\e^{-c_1(L_k)}).
\end{split}
\end{equation}

We say that a bundle $E\to M$ stably splits if
$E\oplus\C^r\simeq\oplus L_j$, where $\C^r$ denotes a trivial
bundle over $M$ and $\oplus L_j$ is a direct sum of complex line
bundles $L_j\to M$. Since $E\simeq\oplus L_k$ (by the Splitting
principle) and $\todd(\C^r)=1$, it follows from \eqref{todd} that
$\todd(E\oplus\C^r)=\todd(E)$. The same is true for $\lclass(E)$
and $\Qy_y(E)$. Now, let $\T$ be a torus and let $E$ be a
$\T$-equivariant vector bundle over $M$. We have the corresponding
$\T$-equivariant versions of \eqref{todd} and
\eqref{ordi-classes},
\begin{equation}\label{todd-lclass}
\begin{split}
&\todd^{\T}(E)=\prod_{k=1}^r\dfrac{c_1^{\T}(L_k)}{1-\e^{-c_1^{\T}(L_k)}}
\quad\mbox{and}\\
&\Qy_y^\T(E)=\todd^{\T}(E)\prod_{k=1}^r\frac{1}{1+y}(1+y
\e^{-c_1^{\T}(L_k)}). \end{split}
\end{equation}


\subsection{Toric manifolds and quantization commutes with reduction}
\label{subse:toricmanifolds} Let $\T=(S^1)^n$ be a torus whose Lie
algebra is denoted by $\liet$. A \textbf{toric manifold} is a
space $(M,\omega,\mu)$, where $(M,\omega)$ is a compact and
connected symplectic $2n$-dimensional manifold on which $\T$ acts
effectively and in a Hamiltonian fashion with moment map $\mu$.
This means that the action of $\T$ on $M$ has trivial kernel, that
the torus action leaves invariant the symplectic form $\omega$ and
that $M$ is provided with a $\T$-equivariant map $\mu:M\to\liet^*$
such that $\imath_{X_{\eta}}\omega=d\langle\mu,\eta\rangle$ for
any $\eta\in\liet$. The function $\langle\mu,\eta\rangle$ is the
component of $\mu$ in the direction of $\eta$ and $X_{\eta}$ is
the vector field on $M$ generated by $\eta$. Examples of toric
manifolds are the complex projective spaces $\CP^n$ and blow ups,
$\CP^n\#\overline{\CP^n}\#\ldots\#\overline{\CP^n}$.\par

Under the moment map, toric manifolds (smooth toric varieties)
correspond to regular polytopes and toric orbifolds (singular
toric varieties) correspond to simple polytopes, modulo
equivariant symplectomorphisms and translations respectively (c.f.
\cite{Del}, \cite{Gin/Gui/Kar}, \cite{Oda}.) The image of a toric
manifold under the moment map is called moment polytope and its
vertices correspond to the fixed point set of the torus action
over the manifold. In general, we say that a symplectic manifold
$(M,\omega)$ is \textbf{prequantizable} if there is a Hermitian
complex line bundle $\lbdle\to M$ such that
$c_1(\lbdle)=[\omega]$. Such a line bundle is a
\textbf{prequantization} of $(M,\omega)$. In the case of a toric
manifold, $M$ is prequantizable if and only if the vertices of its
moment polytope (possibly shifted) are integral lattice points.
(c.f. \cite{Gin/Gui/Kar}, \cite{Gui}.)\par

Let $(M,\omega,\mu)$ be a toric manifold and let $\lbdle\to M$ be
a prequantization of $M$. The theory of geometric quantization
associates a virtual representation $\Q(M)$ of $\T$ to $M$
together with $\lbdle$. The virtual vector space $\Q(M)$ is called
the \textbf{quantization} of the action of $\T$ on $(M,\lbdle)$.
For instance, for the Dolbeault quantization we have
$\Q(M)=\mathrm{kernel}\, D_{\lbdle}\oplus -\mathrm{cokernel}\,
D_{\lbdle}$, where $D_{\lbdle}$ is the Dolbeault operator defined
on the complex $\Omega^*(M)$ of differential forms on $M$ twisted
by $\lbdle$. (c.f. \cite{Gin/Gui/Kar}, \cite{Gui}.) In this case,
the dimension of $\Q(M)$ can be computed by the Atiyah-Singer
index theorem,
\begin{equation}\label{riemann-roch}
\dim\Q(M):= \ind(D_{\lbdle}) = \dint_{M}\e^{c_1(\lbdle)}\todd(M),
\end{equation}
where $\todd(M)$ is the Todd class of the tangent bundle of $M$
and $c_1(\lbdle)$ is the first Chern class of the prequantization
line bundle $\lbdle\to M$.\par

The integer $\dim\Q(M)$ in \eqref{riemann-roch} is equal to the
number of integral lattice points in $\Delta=\mu(M)$ \cite{Dan}.
This is a manifestation of the \textbf{quantization commutes with
reduction} principle in symplectic geometry. To be more precise,
if $\alpha$ is a weight of the torus representation on $\Q(M)$, we
define $\Q(M)^\alpha$ as the subspace of $\Q(M)$ on which $\T$
acts via the character given by $\alpha$. The \emph{quantization
commutes with reduction} principle is stated as
\begin{equation}\label{quant-red}
\dim\Q(M)^\alpha=\dim\Q(M_{\alpha}),
\end{equation}
where $M_{\alpha}$ is the reduced space of $M$ at $\alpha$. The
number $\dim\Q(M)^{\alpha}$ is the \textbf{multiplicity} of the
character given by $\alpha$ in $\Q(M)$ \cite{Gui}. For the
Dolbeault quantization, this multiplicity is 1 or $0$ depending on
whether $\alpha$ is in $\Delta$ or not. Since $M$ is a toric
manifold, the reduced space $M_{\alpha}$ is a point when $\alpha$
is in $\Delta$ and is empty otherwise. Therefore, we have
\begin{equation*}
\dim\Q(M) = \# \mbox{ of integer lattice points in }\Delta =
\sum_{\alpha\in\Delta\cap\Z^n} \dim\Q(M_{\alpha}).
\end{equation*}

\noindent Counting multiplicities boils down to counting integral
lattice points in $\Delta$. We can use Proposition
\ref{ord-decomp} to do this count in terms of an alternating sum
over shifted polarized cones.\par

Now, assume that $\lbdle$ is a $\T$-equivariant prequantization
line bundle for $(M,\omega,\mu)$. This means that the $\T$-action
on $M$ lifts to a $\T$-action on $\lbdle$ by bundle automorphisms
and also that $\lbdle$ is provided with a connection $\nabla$ such
that its equivariant curvature is $\omega-i\mu$. This implies that
$c_1^{\T}(\lbdle)=[\omega-i\mu]$.\footnote{the first
$\T$-equivariant Chern class $c_1^{\T}(\lbdle)$ is of the form
$[\omega+\mu]$ up to a constant. We use the convention
$c_1^{\T}(\lbdle)=[\omega-i\mu]$ for technical reasons.} On the
other hand, by a $\T$-equivariant version of \eqref{riemann-roch}
(see Bismut \cite{Bis}), the character $\chi$ around $0\in\liet$
of this representation of $\T$ on $\Q(M)$ is
\begin{equation}\label{equivriemannroch}
\chi\circ exp = \dint_M\e^{c_1^{\T}(\lbdle)}\todd^{\T}(M),
\end{equation}

\noindent where $exp\,$ denotes the exponential map from $\liet$
to $\T$. Using the Atiyah-Bott Berline-Vergne localization formula
in equivariant cohomology (c.f. \cite{Gin/Gui/Kar}), the right
hand side of \eqref{equivriemannroch} can be expressed as a finite
sum over the fixed points $F$ of the torus action on $M$,
\begin{equation}\label{abbv-riemannroch}
\dint_M\e^{c_1^{\T}(\lbdle)}\todd^{\T}(M)=\dsum_{F}\e^{-i\langle\mu(F),
u\rangle}\prod_{j=1}^n\frac{1}{1-\e^{-i\langle\alpha_{j,F},u\rangle}}.
\end{equation}

\noindent Combining \eqref{abbv-riemannroch} with the notion of
polarization gives us a way of counting the lattice points of the
moment polytope of a toric manifold in terms of the lattice points
of cones based on the vertices of the polytope. This insight
provided by the localization formula in the context of
quantization commutes with reduction for toric manifolds inspires
us to give a weighted polar decomposition for any simple polytope,
not necessarily integral. The proof of this result is completely
independent from the use of the localization formula and is shown
in Section \S \ref{se:polar-decomposition}.\par

\section{Weighted quantization commutes with reduction}
\label{se:wqcr}
\begin{defi}\label{de:chi_y characteristic}\cite{Hir} The equivariant
Hirzebruch $\chi_y$-characteristic of an equivariant complex line
bundle $\lbdle$ over a toric manifold $M$ is
\begin{equation}\label{eq:chi_y characteristic}
\chi_y(M,\lbdle)=\int_M\e^{c_1^\T(\lbdle)}\Qy_y^\T(M).
\end{equation}
\end{defi}

Motivated by \eqref{abbv-riemannroch}, we replace $\todd^\T(M)$ by
$\Qy_y^\T(M)$ in \eqref{equivriemannroch} and use the Atiyah-Bott
Berline-Vergne localization formula again to compute
\eqref{eq:chi_y characteristic}.

Let $(M,\omega,\mu)$ be a dimension $2n$ toric manifold with
moment polytope $\mu(M)\!=\!\Delta$. Assume that the vertices of
$\Delta$ are in $\Z^n$, then $M$ is prequantizable. Let $\lbdle\to
M$ be a prequantization complex line bundle over $M$ and assume
furthermore that $\lbdle$ is $\T$-equivariant. The equivariant
Chern class $c_1^{\T}(\lbdle)$ is equal to $[\omega-i\mu]$. The
complexified tangent bundle of $M$ stably splits into a direct sum
of complex line bundles $\oplus L_j$ which we assume are also
$\T$-equivariant. Then the $\T$-equivariant Hirzebruch class
$\Qy_y^\T(M)$ is
\begin{equation*}
\Qy_y^\T(M)=\todd^{\T}(M)\prod_{j=1}^n\frac{1}{1+y}(1+y
\e^{-c_1^{\T}(L_j)}),
\end{equation*}
\noindent where $y$ is any real number not equal to $-1$. By the
Atiyah-Bott Berline-Vergne localization formula, we get
\begin{align}\label{abbv1}
\begin{split}
\dint_M \e^{c_1^{\T}(\lbdle)}&\Qy_y^\T(M)\\
                  &= \dsum_{F\in M^T}\frac{\e^{c_1^{\T}(\lbdle)}|_F
                  (\todd^{\T}(M)\prod_{j=1}^n \frac{1}{1+y}
                  (1+y\e^{-c_1^{\T}(L_j)}))|_F}{e^{\T}(M)|_F}\\
                  &= \dsum_{F\in M^T} \e^{c_1^{\T}(\lbdle)}|_F \prod_{j=1}^n \left.
                  \left(\frac{1}{1+y}\left(\frac{1+y\e^{-c_1^{\T}(L_j)}}
{1-\e^{-c_1^{\T}(L_j)}}\right)\right)\right|_F.
\end{split}
\end{align}

Let $u$ be a generic element of the Lie algebra $\liet\cong i\R^n$
and $v=\mu(F)$ be the vertex of the moment polytope $\Delta$
corresponding to the fixed point $F$. For $1\le j\le n$, let
$\alpha_{j,v}$ be edge vectors at $v$ normalized as in Subsection
\S\S\ref{subse:simplepolytopes}. This implies that our choice of
such vectors is unique. These edge vectors can be thought of as
the weights of the isotropic representation of $\T$ on the normal
bundle $N_F=T_FM$ at $F$ \cite{Gui}. The restriction of the
equivariant Chern class $c_1^{\T}(L_j)$ to $F$ is
$c_1^{\T}(L_j)|_F=i\langle\alpha_{j,v},u\rangle$ for $1\le j\le
n$. Similarly, the restriction of $c_1^{\T}(\lbdle)$ to $F$ is
$c_1(\lbdle)|_F=-i\langle v,u\rangle$. Then, \eqref{abbv1} can be
written as
\begin{equation}\label{abbv2}
\dint_M \e^{c_1^{\T}(\lbdle)}\Qy_y^\T(M) = \dsum_{v}\e^{-i\langle
v,u\rangle}\prod_{j=1}^n\left(
\frac{1}{1+y}\left(\frac{1+y\e^{-i\langle\alpha_{j,v},u\rangle}}
{1-\e^{-i\langle\alpha_{j,v},u\rangle}}\right)\right)\, ,
\end{equation}

\noindent where the sum is over all vertices of $\Delta$. We pick
$\xi\in\R^n\ssminus\{0\}$ such that
$\langle\alpha_{j,v},\xi\rangle\neq 0$ for $1\le j\le n$ and all
vertices $v$ in $\Delta$. As seen on Subsection
\S\S\ref{subse:simplepolytopes}, we can always find such a $\xi$.
Letting $u=i\xi$, we get
\begin{equation}\label{abbv3}
\begin{split}
\dsum_{v}\e^{-i\langle v,i\xi\rangle}\prod_{j=1}^n&\left(
\frac{1}{1+y}
\left(\frac{1+y\e^{-i\langle\alpha_{j,v},i\xi\rangle}}
{1-\e^{-i\langle\alpha_{j,v},i\xi\rangle}}\right)\right)\\
&=\dsum_{v}\e^{\langle v,\xi\rangle}
\prod_{j=1}^n\left(\left[\frac{1}{1+y}+\frac{y\e^{\langle\alpha_{j,v},
\xi\rangle}}{1+y}\right]\frac{1}{1-\e^{\langle\alpha_{j,v},
\xi\rangle}}\right).
\end{split}
\end{equation}

\noindent Denote $\e^{\langle\alpha_{j,v},\xi\rangle}$ by
$z_{j,v}$. If $\langle\alpha_{j,v},\xi\rangle<0$, we set
$\alpha^\sharp_{j,v}=\alpha_{j,v}$ and have
$z_{j,v}=\e^{\langle\alpha^\sharp_{j,v},\xi\rangle}$, so that
$\vert z_{j,v}\vert<1$ and
$$\frac{1}{1-z_{j,v}}=1+z_{j,v}+z_{j,v}^2+\cdots=\dsum_{m=0}^\infty z_{j,v}^m.$$
If $\langle\alpha_{j,v},\xi\rangle>0$, we set
$\alpha^\sharp_{j,v}=-\alpha_{j,v}$ and have
$z^{-1}_{j,v}=\e^{\langle\alpha^\sharp_{j,v},\xi\rangle}$, so that
$\vert z^{-1}_{jv}\vert<1$ and
$$\frac{1}{1-z_{j,v}}=-z_{j,v}^{-1}\frac{1}{1-z_{j,v}^{-1}}=-z_{j,v}^{-1}
(1+z_{j,v}^{-1}+z_{j,v}^{-2}+\cdots)=-\dsum_{m=1}^\infty
z_{j,v}^{-m}.$$

\noindent Thus, when $\vert z_{j,v}\vert<1$,
\begin{equation}\label{unflipped}
\left(\frac{1}{1+y}+\frac{yz_{j,v}}{1+y}\right)\frac{1}{1-z_{j,v}}=\frac{1}{1+y}+
\dsum_{m=1}^\infty z_{j,v}^m,
\end{equation}
whereas, if $\vert z_{j,v}^{-1}\vert<1$,
\begin{equation}\label{flipped}
\left(\frac{1}{1+y}+\frac{yz_{j,v}}{1+y}\right)\frac{1}{1-z_{j,v}}=
-\left(\frac{y}{1+y}+\dsum_{m=1}^\infty z_{j,v}^{-m}\right).
\end{equation}
\noindent Now, denote $\e^{\langle\alpha^\sharp_{j,v},\xi\rangle}$
by $z_{j,v}$, regardless $\alpha^\sharp_{j,v}=\alpha_{j,v}$ or
$\alpha^\sharp_{j,v}=-\alpha_{j,v}$. Plugging \eqref{unflipped}
and \eqref{flipped} in \eqref{abbv3}, we obtain
\begin{equation}\label{sumprod2}
\dsum_{v}(-1)^{\#v}\e^{\langle v,\xi\rangle}
\prod_{k=1}^{n-\#v}\!\left(\!\frac{1}{1+y}\!+\!\dsum_{m_k=1}^\infty
z_{j_k,v}^{m_k}\right)\prod_{l=1}^{\#v}
\left(\!\frac{y}{1+y}\!+\!\dsum_{m_l=1}^\infty
z_{j_l,v}^{m_l}\right),
\end{equation}

\noindent where $(-1)^{\#v}$ denotes the number of $\alpha_{j,v}$
flipped by the polarizing vector $\xi$. Setting
$v=(v_1,\cdots,v_n)$, $\xi=(\xi_1,\cdots,\xi_n)$ and
$\alpha^\sharp_{j,v}=(\alpha^{\sharp 1}_{j,v},\cdots,
\alpha^{\sharp n}_{j,v})$, we have
\begin{equation*}
\begin{split}
\e^{\langle v,\xi\rangle}&=\e^{v_1\xi_1+\cdots+v_n\xi_n}\quad
\mbox{and}\\
z_{j_k,v}^{m_k}&=\e^{m_k\langle\alpha^\sharp_{j_k,v},\xi\rangle}
                  =\e^{m_k(\alpha^{\sharp 1}_{j,v}\xi_1+\cdots+
\alpha^{\sharp n}_{j,v}\xi_n)}.
\end{split}
\end{equation*}
\noindent Let $z_j=\e^{\xi_j}$. Formula \eqref{sumprod2} is then
equal to
\begin{gather}\label{sumprod3}
  \begin{split}
    \dsum_{v}(-1)^{\#v}z_1^{v_1}\cdots z_n^{v_n}&\prod_{k=1}^{n-\#v}
    \left(\frac{1}{1+y}+\dsum_{m_k=1}^\infty
    z_1^{m_k\alpha^{\sharp 1}_{j_k,v}}\cdots
    z_n^{m_k\alpha^{\sharp n}_{j_k,v}}\right)\\
    & \,\,\,\prod_{l=1}^{\#v}\left(\frac{y}{1+y}+\dsum_{m_l=1}^\infty
    z_1^{m_l\alpha^{\sharp 1}_{j_l,v}}\cdots
    z_n^{m_l\alpha^{\sharp n}_{j_l,v}}\right)
  \end{split}\\
  \begin{align*}
    =\dsum_{v}(-1)^{\#v}z_1^{v_1}\cdots z_n^{v_n}
    \dsum_m w_{v}(m)&z_1^{m_1\alpha^{\sharp 1}_{1,v}+
    \cdots+m_n\alpha^{\sharp 1}_{n,v}}\\
    & z_2^{m_1\alpha^{\sharp 2}_{1,v}+
    \cdots+m_n\alpha^{\sharp 2}_{n,v}}\\
    &\hspace{1.4cm}\vdots\\
    & z_n^{m_1\alpha^{\sharp n}_{1,v}+
    \cdots+m_n\alpha^{\sharp n}_{n,v}},\\
  \end{align*}
\end{gather}
\noindent where $m=(m_1,\cdots,m_n)\in\Z^n_{\ge 0}$ and $w_{v}(m)$
is equal to
\begin{equation}\label{weight-ver1}
w_v(m)=\left(\frac{1}{1+y}\right)^{r_{1,v}(m)}\left(\frac{y}{1+y}\right)^{r_{2,v}(m)}.
\end{equation}

\noindent The exponent $r_{1,v}(m)$ above is the number of $m_j$
in $m$ which correspond to $\alpha_{j,v}^\sharp=\alpha_{j,v}$ and
$r_{2,v}(m)$ is the number of $m_j$ in $m$ that correspond to
$\alpha_{j,v}^\sharp=-\alpha_{j,v}$. We write $z_1^{v_1}\cdots
z_n^{v_n}$ as $z^{v}$ and the long factor in the right hand side
of Equation \eqref{sumprod3} as
$$z_1^{m_1\alpha^{\sharp 1}_{1,v}+
    \cdots+m_n\alpha^{\sharp 1}_{n,v}}z_2^{m_1\alpha^{\sharp 2}_{1,v}+
    \cdots+m_n\alpha^{\sharp 2}_{n,v}}\ldots z_n^{m_1\alpha^{\sharp n}_{1,v}+
    \cdots+m_n\alpha^{\sharp n}_{n,v}} = z^{\sum_j m_j\alpha^\sharp_{j,v}}.$$

\noindent Note that each point $p=(p_1,\ldots,p_n)$ in
$\bfC^\sharp_v\cap\Z^n$ can be written as $p=v+\sum_j
m_j\alpha^\sharp_{j,v}$ with $m_j\in\Z_{\ge 0}$ for $1\le j\le n$,
since $\{\alpha^\sharp_{j,v}\}_{j=1}^n$ form a $\Z$-basis of
$\Z^n$. We write $z^p$ for $z_1^{p_1}\ldots z_n^{p_n}$ and set
$r_{1,v}(p)=r_{1,v}(m)$ and $r_{2,v}(p)=r_{2,v}(m)$ so that
$w_v(p)=w_v(m)$. If $p\notin\bfC^\sharp_v$ we set $w_v(p)=0$. Note
that $r_{1,v}(p)+r_{2,v}(p)$ is equal to the codimension of the
smallest dimensional face in $\bfC^\sharp_v$ containing $p$. With
this notation, we can write \eqref{sumprod3} as
\begin{equation}%
\label{eq:doublesum}
\dsum_{v}(-1)^{\#v}\dsum_{p\in\bfC^\sharp_v\cap\Z^n} w_{v}(p)
z^{p}.
\end{equation}

Let $p\in\Delta\cap\Z^n$ and let $\xi_1$ be a polarizing vector
for $\Delta$. The point $p$ may belong to more than one polarized
tangent cone. Since $p$ is fixed, we have
\begin{equation}%
\label{eq:sum_p_fixed}%
\dsum_{v}(-1)^{\#v} w_{v}(p) z^{p} = \left(\dsum_{v}(-1)^{\#v}
w_{v}(p)\right) z^{p}.
\end{equation}

\noindent To determine the second sum in \eqref{eq:sum_p_fixed},
we define $c(p)$ to be the codimension of the smallest dimensional
face in $\Delta$ containing $p$. If $p$ is in the interior of
$\Delta$ we have $c(p)=0$. Polarizing with respect to $\xi_1$ (in
fact, with respect to any polarizing vector), we have that the
polarized tangent cones of $\Delta$ are all pointing away from $p$
except one which contains $p$ in its interior. Suppose that such
cone is $\bfC^\sharp_v$ for some vertex $v$ of $\Delta$. Then
$r_{1,v}(p) = 0$ and $ r_{2,v}(p) = 0$ since $r_{1,v}(p)+
r_{2,v}(p) = c(p) = 0$. If $c(p)>0$, the point $p$ may belong to
more than one polarized tangent cone. When $c(p)=1$, there are two
cases: the sum $\sum_v (-1)^{\#v} w_v(p)$ is equal to
$\frac{1}{1+y} + \mathit{zeros}$ or equal to $1 - \frac{y}{1+y} +
\mathit{zeros}$. In either case, we have $\sum_v (-1)^{\#v} w_v(p)
= \frac{1}{1+y}$. In general, we can see that (after
simplifications) there are basically three cases: the sum $\sum_v
(-1)^{\#v} w_v(p)$ is equal to $(\frac{1}{1+y})^{c(p)} +
\mathit{zeros}$ or equal to $(1 - \frac{y}{1+y})^{c(p)} +
\mathit{zeros}$ or equal to
$(\frac{1}{1+y})^{c(p)-1}-(\frac{1}{1+y})^{c(p)-1}(\frac{y}{1+y})
+ \mathit{zeros}$. In any case, we have $\sum_v (-1)^{\#v} w_v(p)
= (\frac{1}{1+y})^{c(p)}$. A different polarizing vector $\xi_2$
have the effect of just permuting the type of polarized tangent
cones (that is, cones with no edges flipped, with one edge
flipped, with two edges flipped and so on) over the vertices of
$\Delta$. In conclusion, regardless the polarizing vector $\xi$,
we always obtain
\begin{equation}%
\label{summand}%
\dsum_{v}(-1)^{\#v} w_{v}(p) = \left(\frac{1}{1+y}\right)^{c(p)}.
\end{equation}

From \eqref{sumprod3} on, our arguments apply to any regular
integral polytope in $\R^n$. Therefore, by summing over all the
points of $\Delta\cap\Z^n$, we have proved the following
statement.

\begin{teo}%
\label{th:1}%
Let $\Delta$ be a regular integral polytope in $\R^n$. For any
polarizing vector $\xi\in\R^n$, we have
\begin{equation}\label{sumprod-final}%
\dsum_{v}(-1)^{\#v}\dsum_{p\in\bfC^\sharp_v\cap\Z^n} w_{v}(p)
z^{p} = \dsum_{p\in\Delta\cap\Z^n}
\left(\frac{1}{1+y}\right)^{c(p)}z^p.
\end{equation}
\end{teo}

By considering equations from \eqref{eq:chi_y characteristic} to
\eqref{sumprod-final} together with the work in between, we get
\begin{equation}\label{eq:chi_y weighted}
\chi_y(M,\lbdle) = \dint_M \e^{c_1^{\T}(\lbdle)}\Qy_y^\T(M) =
\dsum_{p\in\Delta\cap\Z^n} \left(\frac{1}{1+y}\right)^{c(p)}z^p.
\end{equation}

The equivariant Hirzebruch $\chi_y$-characteristic of an
equivariant complex line bundle $\lbdle$ over a toric manifold $M$
is in general a polynomial in $\C[z_1^{\pm 1},\ldots,z_n^{\pm
1}]$, as we can see from \eqref{eq:chi_y weighted}. This
equivariant characteristic number is not in general the virtual
index of an elliptic operator on $M$. However, when $y=0$, we have
$\Qy_y^\T(M)=\todd^\T(M)$ and \eqref{eq:chi_y characteristic} is
the equivariant index of the twisted Dolbeault operator $D_\lbdle$
mentioned in Subsection \S\ref{subse:toricmanifolds}. When $y=1$,
we have $\Qy_y^\T(M)=\lclass^\T(M)$ and \eqref{eq:chi_y
characteristic} is in this case the equivariant index of the
twisted signature operator on $M$ \cite{La/Mi}. Because of this,
we can interpret \eqref{eq:chi_y weighted} as a weighted version
of the quantization commutes with reduction principle for a toric
manifold.

\begin{teo}%
\label{th:2}%
Let $M$ be any toric manifold whose moment polytope is denoted by
$\Delta$ and let $\alpha\in\Delta\cap\Z^n$ be any weight lattice
of the corresponding torus action on $M$. Then
\begin{equation}\label{maintheorem2}
\chi_y(M,\lbdle)^\alpha=\left\{%
\begin{array}{cl}
  \left(\frac{1}{1+y}\right)^{c(\alpha)} & \mbox{if }\alpha\in\Delta \\[2ex]
  0 & \mbox{if }\alpha\notin\Delta \\
\end{array}%
\right. .
\end{equation}
\end{teo}

\begin{proof} It follows from \eqref{eq:chi_y
weighted} and the preceding discussion.
\end{proof}




\begin{rem}\label{re:4} In \cite{TZ}, Y. Tian and W. Zhang perform
an analytic approach to the quantization commutes with reduction
principle and get a related weighted multiplicity formula for
Spin$^c$-complexes twisted by certain exterior power bundles.
\end{rem}

\section{Weighted polar decomposition for a simple polytope}
\label{se:polar-decomposition}
Formula \eqref{sumprod-final} is valid for any regular integral
polytope in $\R^n$. In analogy with \eqref{oldformula}, we can
write \eqref{sumprod-final} using characteristic functions and
extend the result to any simple polytope, not necessarily
integral.\par

Let $\Delta$ be a simple polytope in $\R^n$ and let $x$ be any
point in this polytope. Let $c(x)$ be the codimension of the
smallest dimensional face in $\Delta$ containing $x$. We define
the following weighted characteristic function over $\Delta$:
\begin{equation}\label{weightedfunction-polytope}
\mathbf{1}^w_{\Delta}(x)=\left\{%
\begin{array}{cc}
  \left(\frac{1}{1+y}\right)^{c(x)} & \mbox{if }x\in\Delta \\[2ex]
  0 & \mbox{if }x\notin\Delta \\
\end{array}%
\right. .
\end{equation}

\noindent Given a polarizing vector $\xi$ in $\R^n$, let
$\bfC^\sharp_v$ be the polarized cone corresponding to the vertex
$v$. If $x$ is in $\bfC^\sharp_v$, let $c_{v}(x)$ be the
codimension of the smallest dimensional face $F$ in
$\bfC^\sharp_v$ containing $x$. We define $r_{1,v}(x)$ and
$r_{2,v}(x)$ as before; namely, $r_{1,v}(x)$ is the number of
unflipped edges emanating from $v$ which do not belong to $F$, and
$r_{2,v}(x)$ is the number of flipped edges emanating from $v$
which are not in $F$. Therefore, $r_{1,v}(x)+r_{2,v} (x)= c_v(x)$.
For any $y\neq-1$, we set
\begin{equation}\label{multiplicity}
w_{v}(x)=\left(\frac{1}{1+y}\right)^{r_1(x)}\left(\frac{y}{1+y}\right)^{r_2(x)}.
\end{equation}
We define the following weighted characteristic function over
$\bfC^\sharp_v$:
\begin{equation}\label{weightedfunction-cone}
\mathbf{1}^{w_{v}}_{\bfC^\sharp_v}(x)=\left\{%
\begin{array}{cc}
  w_{v}(x) & \mbox{if }x\in \bfC^\sharp_v \\[2ex]
  0 & \mbox{if }x\notin \bfC^\sharp_v \\
\end{array}%
\right. .
\end{equation}

\begin{pro}\label{weig-decomp} For any simple polytope $\Delta$ and for
any choice of polarizing vector, we have
\begin{equation}\label{new-weighted-law}
\mathbf{1}^w_{\Delta}=\sum_{v}(-1)^{\#v}\mathbf{1}^{w_{v}}_{\bfC^\sharp_v},
\end{equation}
where $\mathbf{1}^w_{\Delta}$ and
$\mathbf{1}^{w_{v}}_{\bfC^\sharp_v}$ are the weighted
characteristic functions over $\Delta$ and $\bfC^\sharp_v$ defined
in~(\ref{weightedfunction-polytope})
and~(\ref{weightedfunction-cone}) , where the sum is over the
vertices of $\Delta$, and where $\#v$ denotes the number of edge
vectors at $v$ flipped by the polarizing process
\eqref{polarization}.
\end{pro}


The proof of this result collects the main ideas from the proof of
Theorem \ref{th:1}. Notice that when $y=0$, we obtain the ordinary
polar decomposition stated in Proposition \ref{ord-decomp}. When
$y=1$, the weighted characteristic functions of the polarized
tangent cones are the same at all vertices. In this case, we can
write \eqref{new-weighted-law} as in \cite{KSW2} and \cite{KSW1},
\begin{equation*}\label{old-weighted-law}
\mathbf{1}^w_{\Delta}=\sum_{v}(-1)^{\#v}\mathbf{1}^w_{\bfC^\sharp_v}\,.
\end{equation*}

\begin{figure}[h]
\centering
\includegraphics[scale=.85]{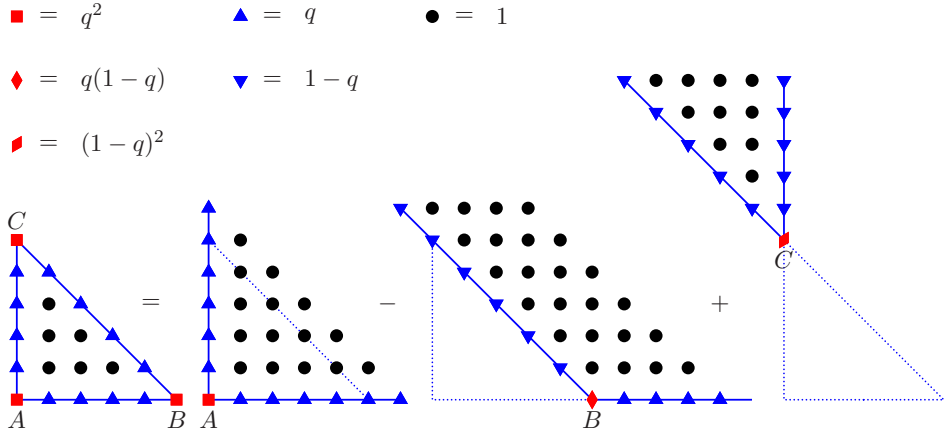}
\caption{Weighted polar decomposition of a lattice
triangle}\label{law-pic}
\end{figure}

We illustrate Proposition~\ref{weig-decomp} for a triangle in
Figure~\ref{law-pic}.

\noindent {\it Proof.} The proof has two steps. First, for each
$x\in\R^n$ we choose a polarizing vector $\xi_x$ such that
\eqref{new-weighted-law} holds. For this, we show that $x$ appears
in at most one polarized tangent cone in the corresponding polar
decomposition of $\Delta$ given by $\xi_x$. Second, we show that
the right hand side of \eqref{new-weighted-law} is independent of
the choice of polarization.\par

Notice that for any polarization there is always a unique cone
with no edges flipped at all. More explicitly, if $x$ is in the
exterior of $\Delta$, we can take the closest vertex of $\Delta$
to $x$ (if there are more than one such a vertex, pick any one of
them). Say we take $v_1$. Then, we set
$\bfC^\sharp_{v_1}=\bfC_{v_1}$. Clearly $x$ is not in
$\bfC^\sharp_{v_1}$. We polarize according to this cone. Namely,
let $\xi_x$ be a polarizing vector in $\R^n$ that makes
$\bfC^\sharp_{v_1}=\bfC_{v_1}$ and let $E_{\xi_x}$ be a
codimension one plane with normal vector parallel to $\xi_x$.
Moving this plane in a parallel way through the polytope, we flip
the edges of the tangent cones whose vertices lie below
$E_{\xi_x}$ (according to $\xi_x$) until we get to the farthest
vertex $v_2$ away from $v_1$ in the direction of $\xi_x$. There is
only one such a vertex since $\xi_x$ is a polarizing vector. The
corresponding polarized tangent cone $\bfC^\sharp_{v_2}$ has all
their edges flipped with respect to $\bfC_{v_2}$ and it does not
contain $x$. In general, a polarizing vector flips an edge if this
edge is shared with a previous polarized tangent cone and leaves
that edge unflipped otherwise. Then, for any vertex other than
$v_1$ and $v_2$, the polarized tangent cones have at least one
edge flipped and they all appear pointing away from $x$. We
conclude that $x$ appears in no polarized tangent cone.\par

If $x$ is in $\Delta$, we can pick any vertex $v$ in the face $F$
whose relative interior contains $x$. Then, we set the polarized
tangent cone $\bfC^\sharp_v$ equal to $\bfC_v$ and choose a
polarizing vector $\xi_x$ that makes $\bfC^\sharp_v=\bfC_v$. We
now repeat the process explained above. Again, all other polarized
tangent cones will point away from $x$, with the cone at the
farthest vertex away from $v$ having all their edges flipped.\par

From here on, we follow \cite{KSW2} closely with some
modifications. Let
$$E_1,\cdots,E_N$$
denote all the different codimension one subspaces of ${\R^n}^*$
that are equal to
$$\alpha_{j,v}^{\bot}=\{\eta\in {\R^n}^*\,|\,\langle\eta,\alpha_{j,v}\rangle=0\}$$
for some $j$ and $v$. If no two edges of $\Delta$ are parallel,
then the number $N$ of such hyperplanes is equal to the number of
edges of $\Delta$. A vector $\xi$ is a \emph{polarizing vector} if
and only if it does not belong to any $E_j$. The \emph{polarized
tangent cones} $\bfC^\sharp_v$ only depend on the connected
component of the complement
$${\R^n}^*\ssminus(E_1\cup\cdots\cup E_N)$$
in which $\xi$ lies. Any two polarizing vectors can be connected
by a path $\xi_t$ in ${\R^n}^*$ which crosses the \emph{walls}
$E_j$ one at a time.\par

The second part of the proof consists of showing that the right
hand side of formula \eqref{new-weighted-law} does not change when
the polarizing vector $\xi_t$ crosses a single wall, $E_k$.\par

As $\xi_t$ crosses the wall $E_k$, the sign of the pairing
$\langle\xi_t,\alpha_{j,v}\rangle$ changes exactly if
$E_k=\alpha_{j,v}^{\bot}$. For each vertex $v$, denote by
$S_{v}(x)$ and $S'_{v}(x)$ its contributions to the right hand
side of formula \eqref{new-weighted-law} before and after $\xi_t$
crosses the wall. The vertices for which these contributions
differ are exactly those vertices that lie on edges $e$ of
$\Delta$ which are perpendicular to $E_k$. They come in pairs
because each edge has two endpoints.\par

Let us concentrate on one such an edge, $e$, with endpoints, say,
$u$ and $v$. Let $\alpha_e$ denote an edge vector at $v$ that
points from $v$ to $u$ along $e$. Suppose that the pairing
$\langle\xi_t,\alpha_e\rangle$ flips its sign from negative to
positive as $\xi_t$ crosses the wall; (otherwise we switch the
roles of $v$ and $u$). The \emph{polarized tangent cones} to
$\Delta$ at $v$, before and after $\xi_t$ crosses the wall, are
$$\bfC^\sharp_v=v+\dsum_{j\in I_e}\R_{\ge 0}\alpha^\sharp_{j,v}+\R_{\ge 0}\alpha_e\quad
\mbox{and}\quad (\bfC^\sharp_v)'=v+\dsum_{j\in I_e}\R_{\ge
0}\alpha^\sharp_{j,v}-\R_{\ge 0}\alpha_e,$$ where
$I_e\subset\{1,\cdots,N\}$ encodes the facets that contain $e$.
(The $\alpha^\sharp_{j,v}$ are the same for the different
$\xi_t$'s because the pairings
$\langle\xi_t,\alpha^\sharp_{j,v}\rangle$ do not flip sign when
$\xi_t$ crosses the wall for $j\in I_e$.) Notice that the cones
$\bfC^\sharp_v$ and $(\bfC^\sharp_v)'$ have a common facet and
their union only depends on the edge $e$ and not on the endpoint
$v$. (This uses the assumption that the polytope $\Delta$ is
simple and follows from the fact that $\alpha_{j,u}\in\R_{\ge
0}\alpha_{j,v}+\R\alpha_e$.)\par

To determine the contributions of $v$ to the right hand side of
~(\ref{new-weighted-law}) before and after $\xi_t$ crosses the
wall, we have two cases to analyze. In the first case, the
smallest dimensional face $F$ containing $x$ also contains the
edge $e$. Then, the contributions are
$$
S_v(x)=\varepsilon\left(\dfrac{1}{1+y}\right)^{r_1(x)}
\left(\dfrac{y}{1+y}\right)^{r_2(x)}\quad\mbox{and}\quad
S'_v(x)=-\varepsilon\cdot 0,
$$
where $\varepsilon\in\{-1,1\}$. Their difference is
$$S_v(x)-S'_v(x)=\varepsilon\left(\dfrac{1}{1+y}\right)^{r_1(x)}
\left(\dfrac{y}{1+y}\right)^{r_2(x)}.$$ The contributions of the
other endpoint, $u$, are
$$
S_u(x)=-\varepsilon\cdot 0 \quad\mbox{and}\quad
S'_u(x)=\varepsilon\left(\dfrac{1}{1+y}\right)^{r_1(x)}
\left(\dfrac{y}{1+y}\right)^{r_2(x)},
$$
and their difference is
$$S_u(x)-S'_u(x)=-\varepsilon\left(\dfrac{1}{1+y}\right)^{r_1(x)}
\left(\dfrac{y}{1+y}\right)^{r_2(x)}.$$ Notice that the numbers
$r_1(x)$ and $r_2(x)$ remain the same for both endpoints because
$F$ contains the edge $e$. Hence, the differences $S_v(x)-S'_v(x)$
and $S_u(x)-S'_u(x)$, for the two endpoints $u$ and $v$ of $e$,
sum to zero.\par

In the second case, the edge $e$ is not contained in $F$. Then,
the contributions are
\begin{equation*}
\begin{split}
S_v(x)=\varepsilon\left(\dfrac{1}{1+y}\right)^{r_1(x)}\left(\dfrac{y}{1+y}\right)^{r_2(x)}
\quad\mbox{and} \\
S'_v(x)=-\varepsilon\left(\dfrac{1}{1+y}\right)^{r_1(x)-1}
\left(\dfrac{y}{1+y}\right)^{r_2(x)+1},
\end{split}
\end{equation*}
and their difference is
$$S_v(x)-S'_v(x)=\varepsilon\left(\dfrac{1}{1+y}\right)^{r_1(x)-1}
\left(\dfrac{y}{1+y}\right)^{r_2(x)}.$$ On the other hand, after
crossing the wall, the smallest dimensional face $F'$ in
$\bfC^\sharp_u$ containing $x$ has codimension one less than the
codimension of $F$ in $\bfC^\sharp_v$, before crossing the wall.
Then, a straightforward calculation shows that the contributions
of the other endpoint, $u$, are
$$
S_u(x)=-\varepsilon\cdot 0 \quad\mbox{and}\quad
S'_u(x)=\varepsilon\left(\dfrac{1}{1+y}\right)^{r_1(x)-1}\left(\dfrac{y}{1+y}\right)^{r_2(x)},
$$
and their difference is
$$S_u(x)-S'_u(x)=-\varepsilon\left(\dfrac{1}{1+y}\right)^{r_1(x)-1}
\left(\dfrac{y}{1+y}\right)^{r_2(x)}.$$ Hence, the differences
$S_v(x)-S'_v(x)$ and $S_u(x)-S'_u(x)$, for the two endpoints $u$
and $v$ of $e$, sum to zero again. $\hfill{\Box}$

\bibliography{referencias}

\end{document}